\documentclass{amsart}
\usepackage{amsfonts}

\newtheorem{theorem}{Theorem}
\theoremstyle{plain}

\newtheorem{corollary}{Corollary}

\newtheorem{proposition}{Proposition}
\newtheorem{remark}{Remark}

\numberwithin{equation}{section}
\input{tcilatex}

\begin{document}
\title[Cassels and Greub-Reinboldt Inequalities]{A Generalisation of the
Cassels and Greub-Reinboldt Inequalities in Inner Product Spaces}
\author{Sever S. Dragomir}
\address{School of Computer Science and Mathematics\\
Victoria University of Technology\\
PO Box 14428, MCMC 8001\\
Victoria, Australia.}
\email{sever@matilda.vu.edu.au}
\urladdr{http://rgmia.vu.edu.au/SSDragomirWeb.html}

\begin{abstract}
A generalisation of the Cassels and Greub-Reinboldt inequalities in complex
or real inner product spaces and applications for isotonic \ linear
functionals, integrals and sequences are provided.
\end{abstract}

\maketitle

\section{Introduction}

The following result was proved by J.W.S. Cassels in 1951 (see Appendix 1 of 
\cite{8b}).

\begin{theorem}
\label{t1.1}Let $\overline{\mathbf{a}}=\left( a_{1},\dots ,a_{n}\right) ,$ $%
\overline{\mathbf{b}}=\left( b_{1},\dots ,b_{n}\right) $ be sequences of
positive real numbers and $\overline{\mathbf{w}}=\left( w_{1},\dots
,w_{n}\right) $ a sequence of nonnegative real numbers. Suppose that%
\begin{equation}
m=\min_{i=\overline{1,n}}\left\{ \frac{a_{i}}{b_{i}}\right\} \text{ \ and \ }%
M=\max_{i=\overline{1,n}}\left\{ \frac{a_{i}}{b_{i}}\right\} .  \label{1.1}
\end{equation}%
Then one has the inequality%
\begin{equation}
\frac{\sum_{i=1}^{n}w_{i}a_{i}^{2}\sum_{i=1}^{n}w_{i}b_{i}^{2}}{\left(
\sum_{i=1}^{n}w_{i}a_{i}b_{i}\right) ^{2}}\leq \frac{\left( m+M\right) ^{2}}{%
4mM}.  \label{1.2}
\end{equation}%
The equality holds in (\ref{1.2}) when $w_{1}=\frac{1}{a_{1}b_{1}},$ $w_{n}=%
\frac{1}{a_{n}b_{n}},$ $w_{2}=\cdots =w_{n-1}=0,$ $m=\frac{a_{n}}{b_{1}}$
and $M=\frac{a_{1}}{b_{n}}.$
\end{theorem}

If one assumes that $0<a\leq a_{i}\leq A<\infty $ and $0<b\leq b_{i}\leq
B<\infty $ for each $i\in \left\{ 1,\dots ,n\right\} ,$ then by (\ref{2.2})
we may obtain Greub-Reinboldt's inequality \cite{2b}%
\begin{equation}
\frac{\sum_{i=1}^{n}w_{i}a_{i}^{2}\sum_{i=1}^{n}w_{i}b_{i}^{2}}{\left(
\sum_{i=1}^{n}w_{i}a_{i}b_{i}\right) ^{2}}\leq \frac{\left( ab+AB\right) ^{2}%
}{4abAB}.  \label{1.3}
\end{equation}%
The following \textquotedblleft unweighted\textquotedblright\ Cassels'
inequality also holds%
\begin{equation}
\frac{\sum_{i=1}^{n}a_{i}^{2}\sum_{i=1}^{n}b_{i}^{2}}{\left(
\sum_{i=1}^{n}a_{i}b_{i}\right) ^{2}}\leq \frac{\left( m+M\right) ^{2}}{4mM},
\label{1.4}
\end{equation}%
provided $\overline{\mathbf{a}}$ and $\overline{\mathbf{b}}$ satisfy (\ref%
{1.1}). This inequality will produce the well known P\'{o}lya-Szeg\"{o}
inequality \cite[pp. 57, 213-114]{5b}, \cite[pp. 71-72, 253-255]{3b}:%
\begin{equation}
\frac{\sum_{i=1}^{n}a_{i}^{2}\sum_{i=1}^{n}b_{i}^{2}}{\left(
\sum_{i=1}^{n}a_{i}b_{i}\right) ^{2}}\leq \frac{\left( ab+AB\right) ^{2}}{%
4abAB},  \label{1.5}
\end{equation}%
provided $0<a\leq a_{i}\leq A<\infty $ and $0<b\leq b_{i}\leq B<\infty $ for
each $i\in \left\{ 1,\dots ,n\right\} .$

In \cite{4ba}, C.P. Niculescu proved, amongst others, the following
generalisation of Cassels' inequality:

\begin{theorem}
\label{t1.2}Let $E$ be a vector space endowed with a Hermitian product $%
\left\langle \cdot ,\cdot \right\rangle .$ Then%
\begin{equation}
\frac{\func{Re}\left\langle x,y\right\rangle }{\left\langle x,x\right\rangle
^{\frac{1}{2}}\left\langle y,y\right\rangle ^{\frac{1}{2}}}\geq \frac{2}{%
\sqrt{\frac{\omega }{\Omega }}+\sqrt{\frac{\Omega }{\omega }}}  \label{1.6}
\end{equation}%
for every $x,y\in E$ and every $\omega ,\Omega >0$ for which $\func{Re}%
\left\langle x-\omega y,x-\Omega y\right\rangle \leq 0.$
\end{theorem}

For other reverses of the Cauchy-Bunyakovsky-Schwarz, see the references 
\cite{1b}-\cite{8b}.

In this paper we obtain a generalisation of (\ref{1.6}) for complex numbers $%
\omega $ and $\Omega $ for which $\func{Re}\left( \overline{\omega }\Omega
\right) >0.$ Applications for isotonic linear functionals, integrals and
sequences are also given.

\section{The Results}

The following reverse of Schwarz's inequality in inner product spaces holds.

\begin{theorem}
\label{t2.1}Let $a,A\in \mathbb{K}$ $\left( \mathbb{K}=\mathbb{C},\mathbb{R}%
\right) $ so that $\func{Re}\left( \overline{a}A\right) >0.$ If $x,y\in H$
are such that%
\begin{equation}
\func{Re}\left\langle Ay-x,x-ay\right\rangle \geq 0,  \label{2.1}
\end{equation}%
then one has the inequality%
\begin{equation}
\left\Vert x\right\Vert \left\Vert y\right\Vert \leq \frac{1}{2}\cdot \frac{%
\func{Re}\left[ A\overline{\left\langle x,y\right\rangle }+\overline{a}%
\left\langle x,y\right\rangle \right] }{\left[ \func{Re}\left( \overline{a}%
A\right) \right] ^{\frac{1}{2}}}\leq \frac{1}{2}\cdot \frac{\left\vert
A\right\vert +\left\vert a\right\vert }{\left[ \func{Re}\left( \overline{a}%
A\right) \right] ^{\frac{1}{2}}}\left\vert \left\langle x,y\right\rangle
\right\vert .  \label{2.2}
\end{equation}%
The constant $\frac{1}{2}$ is sharp in both inequalities.
\end{theorem}

\begin{proof}
We have, obviously, that%
\begin{equation*}
I:=\func{Re}\left\langle Ay-x,x-ay\right\rangle =\func{Re}\left[ A\overline{%
\left\langle x,y\right\rangle }+\overline{a}\left\langle x,y\right\rangle %
\right] -\left\Vert x\right\Vert ^{2}-\left[ \func{Re}\left( \overline{a}%
A\right) \right] \left\Vert y\right\Vert ^{2}
\end{equation*}%
and, thus, by (\ref{2.1}), one has%
\begin{equation*}
\left\Vert x\right\Vert ^{2}+\left[ \func{Re}\left( \overline{a}A\right) %
\right] \cdot \left\Vert y\right\Vert ^{2}\leq \func{Re}\left[ A\overline{%
\left\langle x,y\right\rangle }+\overline{a}\left\langle x,y\right\rangle %
\right] ,
\end{equation*}%
giving%
\begin{equation}
\frac{1}{\left[ \func{Re}\left( \overline{a}A\right) \right] ^{\frac{1}{2}}}%
\left\Vert x\right\Vert ^{2}+\left[ \func{Re}\left( \overline{a}A\right) %
\right] ^{\frac{1}{2}}\left\Vert y\right\Vert ^{2}\leq \frac{\func{Re}\left[
A\overline{\left\langle x,y\right\rangle }+\overline{a}\left\langle
x,y\right\rangle \right] }{\left[ \func{Re}\left( \overline{a}A\right) %
\right] ^{\frac{1}{2}}}.  \label{2.3}
\end{equation}%
On the other hand, by the elementary inequality%
\begin{equation*}
\alpha p^{2}+\frac{1}{\alpha }q^{2}\geq 2pq
\end{equation*}%
holding for $p,q\geq 0$ and $\alpha >0,$ we deduce%
\begin{equation}
2\left\Vert x\right\Vert \left\Vert y\right\Vert \leq \frac{1}{\left[ \func{%
Re}\left( \overline{a}A\right) \right] ^{\frac{1}{2}}}\left\Vert
x\right\Vert ^{2}+\left[ \func{Re}\left( \overline{a}A\right) \right] ^{%
\frac{1}{2}}\left\Vert y\right\Vert ^{2}.  \label{2.4}
\end{equation}%
Utilizing (\ref{2.3}) and (\ref{2.4}) we deduce the first part of (\ref{2.2}%
).

The last part is obvious by the fact that for $z\in \mathbb{C},$ $\left\vert 
\func{Re}\left( z\right) \right\vert \leq \left\vert z\right\vert .$

Now, assume that the first inequality in (\ref{2.2}) holds with a constant $%
c>0,$ i.e., 
\begin{equation}
\left\Vert x\right\Vert \left\Vert y\right\Vert \leq c\frac{\func{Re}\left[ A%
\overline{\left\langle x,y\right\rangle }+\overline{a}\left\langle
x,y\right\rangle \right] }{\left[ \func{Re}\left( \overline{a}A\right) %
\right] ^{\frac{1}{2}}},  \label{2.5}
\end{equation}%
where $a,A,x$ and $y$ satisfy (\ref{2.2}).

If we choose $a=A=1,$ $y=x\neq 0,$ then obviously (\ref{2.2}) holds and from
(\ref{2.5}) we may obtain%
\begin{equation*}
\left\Vert x\right\Vert ^{2}\leq 2c\left\Vert x\right\Vert ^{2},
\end{equation*}%
giving $c\geq \frac{1}{2}.$

The theorem is completely proved.
\end{proof}

The following corollary is a natural consequence of the above theorem.

\begin{corollary}
\label{c2.2}Let $m,M>0.$ If $x,y\in H$ are such that%
\begin{equation}
\func{Re}\left\langle My-x,x-my\right\rangle \geq 0,  \label{2.6}
\end{equation}%
then one has the inequality%
\begin{equation}
\left\Vert x\right\Vert \left\Vert y\right\Vert \leq \frac{1}{2}\cdot \frac{%
M+m}{\sqrt{mM}}\func{Re}\left\langle x,y\right\rangle \leq \frac{1}{2}\cdot 
\frac{M+m}{\sqrt{mM}}\left\vert \left\langle x,y\right\rangle \right\vert .
\label{2.7}
\end{equation}%
The constant $\frac{1}{2}$ is sharp in (\ref{2.7}).
\end{corollary}

\begin{remark}
\label{r2.2.1}The inequality (\ref{2.7}) is equivalent to Niculescu's
inequality (\ref{1.6}).
\end{remark}

The following corollary is also obvious.

\begin{corollary}
\label{c2.3}With the assumptions of Corollary \ref{c2.2}, we have%
\begin{align}
0& \leq \left\Vert x\right\Vert \left\Vert y\right\Vert -\left\vert
\left\langle x,y\right\rangle \right\vert \leq \left\Vert x\right\Vert
\left\Vert y\right\Vert -\func{Re}\left\langle x,y\right\rangle  \label{2.8}
\\
& \leq \frac{\left( \sqrt{M}-\sqrt{m}\right) ^{2}}{2\sqrt{mM}}\func{Re}%
\left\langle x,y\right\rangle \leq \frac{\left( \sqrt{M}-\sqrt{m}\right) ^{2}%
}{2\sqrt{mM}}\left\vert \left\langle x,y\right\rangle \right\vert  \notag
\end{align}%
and%
\begin{align}
0& \leq \left\Vert x\right\Vert ^{2}\left\Vert y\right\Vert ^{2}-\left\vert
\left\langle x,y\right\rangle \right\vert ^{2}\leq \left\Vert x\right\Vert
^{2}\left\Vert y\right\Vert ^{2}-\left[ \func{Re}\left\langle
x,y\right\rangle \right] ^{2}  \label{2.9} \\
& \leq \frac{\left( M-m\right) ^{2}}{4mM}\left[ \func{Re}\left\langle
x,y\right\rangle \right] ^{2}\leq \frac{\left( M-m\right) ^{2}}{4mM}%
\left\vert \left\langle x,y\right\rangle \right\vert ^{2}.  \notag
\end{align}
\end{corollary}

\begin{proof}
If we subtract $\func{Re}\left\langle x,y\right\rangle \geq 0$ from the
first inequality in (\ref{2.7}), we get%
\begin{align*}
\left\Vert x\right\Vert \left\Vert y\right\Vert -\func{Re}\left\langle
x,y\right\rangle & \leq \left( \frac{1}{2}\cdot \frac{M+m}{\sqrt{mM}}%
-1\right) \func{Re}\left\langle x,y\right\rangle \\
& =\frac{\left( \sqrt{M}-\sqrt{m}\right) ^{2}}{2\sqrt{mM}}\func{Re}%
\left\langle x,y\right\rangle
\end{align*}%
which proves the third inequality in (\ref{2.8}). The other ones are obvious.

Now, if we square the first inequality in (\ref{2.7}) and then subtract $%
\left[ \func{Re}\left\langle x,y\right\rangle \right] ^{2},$ we get%
\begin{align*}
\left\Vert x\right\Vert ^{2}\left\Vert y\right\Vert ^{2}-\left[ \func{Re}%
\left\langle x,y\right\rangle \right] ^{2}& \leq \left[ \frac{\left(
M+m\right) ^{2}}{4mM}-1\right] \left[ \func{Re}\left\langle x,y\right\rangle %
\right] ^{2} \\
& =\frac{\left( M-m\right) ^{2}}{4mM}\left[ \func{Re}\left\langle
x,y\right\rangle \right] ^{2}
\end{align*}%
which proves the third inequality in (\ref{2.5}). The other ones are obvious.
\end{proof}

\section{Applications for Isotonic Linear Functionals}

Let $F\left( T\right) $ be an algebra of real functions defined on $T$ and $%
L $ a subclass of $F\left( T\right) $ satisfying the conditions:

\begin{enumerate}
\item[(i)] $f,g\in L$ implies $f+g\in L;$

\item[(ii)] $f\in L,$ $\in \mathbb{R}$ implies $\alpha f\in L.$
\end{enumerate}

A functional $A$ defined on $L$ is an \textit{isotonic linear functional} on 
$L$ provided that

\begin{enumerate}
\item[(a)] $A\left( \alpha f+\beta g\right) =\alpha A\left( f\right) +\beta
A\left( g\right) $ for all $\alpha ,\beta \in \mathbb{R}$ and $f,g\in L;$

\item[(aa)] $f\geq g,$ that is, $f\left( t\right) \geq g\left( t\right) $
for all $t\in T,$ implies $A\left( f\right) \geq A\left( g\right) .$
\end{enumerate}

The functional $A$ is \textit{normalised} on $L,$ provided that $\mathbf{1}%
\in L,$ i.e., $\mathbf{1}\left( t\right) =1$ for all $t\in T,$ implies $%
A\left( \mathbf{1}\right) =1.$

Usual examples of isotonic linear functionals are integrals, sums, etc.

Now, suppose that $h\in F\left( T\right) ,$ $h\geq 0$ is given and satisfies
the properties that $fgh\in L,$ $fh\in L,$ $gh\in L$ for all $f,g\in L.$ For
a given isotonic linear functional $A:L\rightarrow \mathbb{R}$ with $A\left(
h\right) >0,$ define the mapping $\left( \cdot ,\cdot \right) _{A,h}:L\times
L\rightarrow \mathbb{R}$ by%
\begin{equation*}
\left( f,g\right) _{A,h}:=\frac{A\left( fgh\right) }{A\left( h\right) }.
\end{equation*}

This functional satisfies the following properties:

\begin{enumerate}
\item[(s)] $\left( f,f\right) _{A,h}\geq 0$ for all $f\in L;$

\item[(ss)] $\left( \alpha f+\beta g,k\right) _{A,h}=\alpha \left(
f,k\right) _{A,h}+\beta \left( g,k\right) _{A,h}$ for all $f,g,k\in L$ and $%
\alpha ,\beta \in \mathbb{R}$;

\item[(sss)] $\left( f,g\right) _{A,h}=\left( g,f\right) _{A,h}$ for all $%
f,g\in L.$
\end{enumerate}

The following proposition holds.

\begin{proposition}
\label{p3.1}Let $f,g,h\in F\left( T\right) $ be such that $fgh\in L,$ $%
f^{2}h\in L,$ $g^{2}h\in L.$ If $m,M>0$ are such that%
\begin{equation}
mg\leq f\leq Mg\text{ on }F\left( T\right) ,  \label{3.2}
\end{equation}%
then for any isotonic linear functional $A:L\rightarrow \mathbb{R}$ with $%
A\left( h\right) >0,$ we have the inequality%
\begin{equation}
1\leq \frac{A\left( f^{2}h\right) A\left( g^{2}h\right) }{A^{2}\left(
fgh\right) }\leq \frac{\left( M+m\right) ^{2}}{4mM}.  \label{3.3}
\end{equation}%
The constant $\frac{1}{4}$ in (\ref{3.3}) is sharp.
\end{proposition}

\begin{proof}
We observe that%
\begin{equation*}
\left( Mg-f,f-mg\right) _{A,h}=A\left[ h\left( Mg-f\right) \left(
f-mg\right) \right] \geq 0.
\end{equation*}%
Applying Corollary \ref{c2.2} for $\left( \cdot ,\cdot \right) _{A,L}$ we get%
\begin{equation*}
1\leq \frac{\left( f,f\right) _{A,h}\left( g,g\right) _{A,h}}{\left(
f,g\right) _{A,h}^{2}}\leq \frac{\left( M+m\right) ^{2}}{4mM},
\end{equation*}%
which is clearly equivalent to (\ref{3.3}).
\end{proof}

The following additive versions of (\ref{3.3}) also hold.

\begin{corollary}
\label{c3.2}With the assumption in Proposition \ref{p3.1}, one has%
\begin{align}
0& \leq \left[ A\left( f^{2}h\right) A\left( g^{2}h\right) \right] ^{\frac{1%
}{2}}-A\left( hfg\right)  \label{3.4} \\
& \leq \frac{\left( \sqrt{M}-\sqrt{m}\right) ^{2}}{2\sqrt{mM}}A\left(
hfg\right)  \notag
\end{align}%
and 
\begin{align}
0& \leq A\left( f^{2}h\right) A\left( g^{2}h\right) -A^{2}\left( fgh\right)
\label{3.5} \\
& \leq \frac{\left( M-m\right) ^{2}}{4mM}A^{2}\left( fgh\right) .
\end{align}
\end{corollary}

\begin{remark}
\label{r3.3}The condition (\ref{3.2}) may be replaced with the weaker
assumption%
\begin{equation}
\left( Mg-f,f-mg\right) _{A,h}\geq 0.  \label{3.6}
\end{equation}
\end{remark}

\begin{remark}
\label{r3.4}With the assumption (\ref{3.2}) or (\ref{3.6}) and if $f,g\in
F\left( T\right) $ with $fg,f^{2},g^{2}\in L,$ then one has the inequality%
\begin{equation}
1\leq \frac{A\left( f^{2}\right) A\left( g^{2}\right) }{A^{2}\left(
fg\right) }\leq \frac{\left( M+m\right) ^{2}}{4mM},  \label{3.7}
\end{equation}%
\begin{align}
0& \leq \left[ A\left( f^{2}\right) A\left( g^{2}\right) \right] ^{\frac{1}{2%
}}-A\left( fg\right)  \label{3.8} \\
& \leq \frac{\left( \sqrt{M}-\sqrt{m}\right) ^{2}}{2\sqrt{mM}}A\left(
fg\right)  \notag
\end{align}%
and%
\begin{equation}
0\leq A\left( f^{2}\right) A\left( g^{2}\right) -A^{2}\left( fg\right) \leq 
\frac{\left( M-m\right) ^{2}}{4mM}A^{2}\left( fg\right) .  \label{3.9}
\end{equation}
\end{remark}

\section{Applications for Integrals}

Let $\left( \Omega ,\Sigma ,\mu \right) $ be a measure space consisting of a
set $\Omega ,$ a $\sigma -$algebra $\Sigma $ of subsets of $\Omega $ and a
countably additive and positive measure on $\Sigma $ with values in $\mathbb{%
R}\cup \left\{ \infty \right\} .$

Denote by $L_{\rho }^{2}\left( \Omega ,\mathbb{K}\right) $ the Hilbert space
of all $\mathbb{K}$-valued functions $f$ defined on $\Omega $ that are $%
2-\rho -$integrable on $\Omega ,$ i.e., $\int_{\Omega }\rho \left( t\right)
\left\vert f\left( s\right) \right\vert ^{2}d\mu \left( s\right) <\infty ,$
where $\rho :\Omega \rightarrow \lbrack 0,\infty )$ is a measurable function
on $\Omega .$

The following proposition contains a counterpart of the weighted
Cauchy-Bunyakovsky-Schwarz's integral inequality.

\begin{proposition}
\label{p4.1}Let $A,a\in \mathbb{K}$ $\left( \mathbb{K}=\mathbb{C},\mathbb{R}%
\right) $ with $\func{Re}\left( \overline{a}A\right) >0$ and $f,g\in L_{\rho
}^{2}\left( \Omega ,\mathbb{K}\right) .$ If 
\begin{equation}
\int_{\Omega }\func{Re}\left[ \left( Ag\left( s\right) -f\left( s\right)
\right) \left( \overline{f\left( s\right) }-\overline{a}\text{ }\overline{g}%
\left( s\right) \right) \right] \rho \left( s\right) d\mu \left( s\right)
\geq 0,  \label{4.1}
\end{equation}%
then one has the inequality%
\begin{align}
& \left[ \int_{\Omega }\rho \left( s\right) \left\vert f\left( s\right)
\right\vert ^{2}d\mu \left( s\right) \int_{\Omega }\rho \left( s\right)
\left\vert g\left( s\right) \right\vert ^{2}d\mu \left( s\right) \right] ^{%
\frac{1}{2}}  \label{4.2} \\
& \leq \frac{1}{2}\cdot \frac{\int_{\Omega }\rho \left( s\right) \func{Re}%
\left[ A\overline{f\left( s\right) }g\left( s\right) +\overline{a}f\left(
s\right) \overline{g\left( s\right) }\right] d\mu \left( s\right) }{\left[ 
\func{Re}\left( \overline{a}A\right) \right] ^{\frac{1}{2}}}  \notag \\
& \leq \frac{1}{2}\cdot \frac{\left\vert A\right\vert +\left\vert
a\right\vert }{\left[ \func{Re}\left( \overline{a}A\right) \right] ^{\frac{1%
}{2}}}\left\vert \int_{\Omega }\rho \left( s\right) f\left( s\right) 
\overline{g\left( s\right) }d\mu \left( s\right) \right\vert .  \notag
\end{align}%
The constant $\frac{1}{2}$ is sharp in (\ref{4.2}).
\end{proposition}

\begin{proof}
Follows by Theorem \ref{t2.1} applied for the inner product $\left\langle
\cdot ,\cdot \right\rangle _{\rho }:=L_{\rho }^{2}\left( \Omega ,\mathbb{K}%
\right) \times L_{\rho }^{2}\left( \Omega ,\mathbb{K}\right) \rightarrow 
\mathbb{K}$,%
\begin{equation*}
\left\langle f,g\right\rangle :=\int_{\Omega }\rho \left( s\right) f\left(
s\right) \overline{g\left( s\right) }d\mu \left( s\right) .
\end{equation*}
\end{proof}

\begin{remark}
\label{r4.2}A sufficient condition for (\ref{4.1}) to hold is 
\begin{equation}
\func{Re}\left[ \left( Ag\left( s\right) -f\left( s\right) \right) \left( 
\overline{f\left( s\right) }-\overline{a}\overline{g\left( s\right) }\right) %
\right] \geq 0\text{ \ for }\mu \text{-a.e. \ }s\in \Omega .  \label{4.3}
\end{equation}
\end{remark}

In the particular case $\rho =1,$ we have the following result.

\begin{corollary}
\label{c4.3}Let $a,A\in \mathbb{K}$ $\left( \mathbb{K}=\mathbb{C},\mathbb{R}%
\right) $ with $\func{Re}\left( \overline{a}A\right) >0$ and $f,g\in
L^{2}\left( \Omega ,\mathbb{K}\right) .$ If%
\begin{equation}
\int_{\Omega }\func{Re}\left[ \left( Ag\left( s\right) -f\left( s\right)
\right) \left( \overline{f\left( s\right) }-\overline{a}\overline{g\left(
s\right) }\right) \right] d\mu \left( s\right) \geq 0,  \label{4.4}
\end{equation}%
then one has the inequality%
\begin{align}
& \left[ \int_{\Omega }\left\vert f\left( s\right) \right\vert ^{2}d\mu
\left( s\right) \int_{\Omega }\left\vert g\left( s\right) \right\vert
^{2}d\mu \left( s\right) \right] ^{\frac{1}{2}}  \label{4.5} \\
& \leq \frac{1}{2}\cdot \frac{\int_{\Omega }\func{Re}\left[ A\overline{%
f\left( s\right) }g\left( s\right) +\overline{a}f\left( s\right) \overline{%
g\left( s\right) }\right] d\mu \left( s\right) }{\left[ \func{Re}\left( 
\overline{a}A\right) \right] ^{\frac{1}{2}}}  \notag \\
& \leq \frac{1}{2}\cdot \frac{\left\vert A\right\vert +\left\vert
a\right\vert }{\left[ \func{Re}\left( \overline{a}A\right) \right] ^{\frac{1%
}{2}}}\left\vert \int_{\Omega }f\left( s\right) \overline{g\left( s\right) }%
d\mu \left( s\right) \right\vert .  \notag
\end{align}
\end{corollary}

\begin{remark}
\label{r4.4}If $\mathbb{K}=\mathbb{R},$ then a sufficient condition for
either (\ref{4.1}) or (\ref{4.4}) \ to hold is 
\begin{equation}
ag\left( s\right) \leq f\left( s\right) \leq Ag\left( s\right) \text{ \ for
\ }\mu \text{-a.e. }s\in \Omega ,  \label{4.6}
\end{equation}%
where, in this case, $a,A\in \mathbb{R}$ with $A>a>0.$
\end{remark}

When $a,A$ are real positive constants, then the following proposition holds.

\begin{proposition}
\label{p4.5}Let $m,M>0.$ If $f,g\in L_{\rho }^{2}\left( \Omega ,\mathbb{K}%
\right) $ such that%
\begin{equation}
\int_{\Omega }\rho \left( s\right) \func{Re}\left[ \left( Mg\left( s\right)
-f\left( s\right) \right) \left( \overline{f\left( s\right) }-m\overline{g}%
\left( s\right) \right) \right] d\mu \left( s\right) \geq 0  \label{4.7}
\end{equation}%
then one has the inequality%
\begin{multline}
\left[ \int_{\Omega }\rho \left( s\right) \left\vert f\left( s\right)
\right\vert ^{2}d\mu \left( s\right) \int_{\Omega }\rho \left( s\right)
\left\vert g\left( s\right) \right\vert ^{2}d\mu \left( s\right) \right] ^{%
\frac{1}{2}}  \label{4.8} \\
\leq \frac{1}{2}\cdot \frac{M+m}{\sqrt{mM}}\int_{\Omega }\rho \left(
s\right) \func{Re}\left[ f\left( s\right) \overline{g\left( s\right) }\right]
d\mu \left( s\right) .
\end{multline}
\end{proposition}

The proof follows by Corollary \ref{c2.2} applied for the inner product 
\begin{equation*}
\left\langle f,g\right\rangle _{\rho }:=\int_{\Omega }\rho \left( s\right)
f\left( s\right) \overline{g\left( s\right) }d\mu \left( s\right) .
\end{equation*}

The following additive versions also hold.

\begin{corollary}
\label{c4.6}With the assumptions in Proposition \ref{p4.5}, one has the
inequalities%
\begin{align}
0& \leq \left[ \int_{\Omega }\rho \left( s\right) \left\vert f\left(
s\right) \right\vert ^{2}d\mu \left( s\right) \int_{\Omega }\rho \left(
s\right) \left\vert g\left( s\right) \right\vert ^{2}d\mu \left( s\right) %
\right] ^{\frac{1}{2}}  \label{4.9} \\
& \ \ \ \ \ \ \ \ \ \ \ \ \ \ \ \ \ \ \ \ \ \ \ \ \ -\int_{\Omega }\rho
\left( s\right) \func{Re}\left[ f\left( s\right) \overline{g\left( s\right) }%
\right] d\mu \left( s\right)  \notag \\
& \leq \frac{\left( \sqrt{M}-\sqrt{m}\right) ^{2}}{2\sqrt{mM}}\int_{\Omega
}\rho \left( s\right) \func{Re}\left[ f\left( s\right) \overline{g\left(
s\right) }\right] d\mu \left( s\right)  \notag
\end{align}%
and 
\begin{align}
0& \leq \int_{\Omega }\rho \left( s\right) \left\vert f\left( s\right)
\right\vert ^{2}d\mu \left( s\right) \int_{\Omega }\rho \left( s\right)
\left\vert g\left( s\right) \right\vert ^{2}d\mu \left( s\right)
\label{4.10} \\
& \ \ \ \ \ \ \ \ \ \ \ \ \ \ \ \ \ \ \ \ \ \ \ \ -\left( \int_{\Omega }\rho
\left( s\right) \func{Re}\left[ f\left( s\right) \overline{g\left( s\right) }%
\right] d\mu \left( s\right) \right) ^{2}  \notag \\
& \leq \frac{\left( M-m\right) ^{2}}{4mM}\left( \int_{\Omega }\rho \left(
s\right) \func{Re}\left[ f\left( s\right) \overline{g\left( s\right) }\right]
d\mu \left( s\right) \right) ^{2}.  \notag
\end{align}
\end{corollary}

\begin{remark}
\label{r4.7}If $\mathbb{K}=\mathbb{R}$, a sufficient condition for (\ref{4.1}%
) to hold is 
\begin{equation}
mg\left( s\right) \leq f\left( s\right) \leq Mg\left( s\right) \text{ \ for
\ }\mu \text{-a.e. }s\in \Omega ,  \label{4.11}
\end{equation}%
where $M>m>0.$
\end{remark}

\section{Applications for Sequences}

For a given sequence $\left( w_{i}\right) _{i\in \mathbb{N}}$ of nonnegative
real numbers, consider the Hilbert space $\ell _{w}^{2}\left( \mathbb{K}%
\right) ,$ $\left( \mathbb{K}=\mathbb{C},\mathbb{R}\right) ,$ where 
\begin{equation}
\ell _{w}^{2}\left( \mathbb{K}\right) :=\left\{ \overline{x}=\left(
x_{i}\right) _{i\in \mathbb{N}}\subset \mathbb{K}\left\vert
\sum_{i=0}^{\infty }w_{i}\left\vert x_{i}\right\vert ^{2}<\infty \right.
\right\} .  \label{5.1}
\end{equation}

The following proposition that provides a counterpart of the weighted
Cauchy-Bunyakovsky-Schwarz inequality for complex numbers holds.

\begin{proposition}
\label{p5.1}Let $a,A\in \mathbb{K}$ with $\func{Re}\left( \overline{a}%
A\right) >0$ and $\overline{x},\overline{y}\in \ell _{w}^{2}\left( \mathbb{K}%
\right) .$ If%
\begin{equation}
\sum_{i=0}^{\infty }w_{i}\func{Re}\left[ \left( Ay_{i}-x_{i}\right) \left( 
\overline{x_{i}}-\overline{a}\overline{y_{i}}\right) \right] \geq 0,
\label{5.2}
\end{equation}%
then one has the inequality%
\begin{align}
\left[ \sum_{i=0}^{\infty }w_{i}\left\vert x_{i}\right\vert
^{2}\sum_{i=0}^{\infty }w_{i}\left\vert y_{i}\right\vert ^{2}\right] ^{\frac{%
1}{2}}& \leq \frac{1}{2}\cdot \frac{\sum_{i=0}^{\infty }w_{i}\func{Re}\left[
A\overline{x_{i}}y_{i}+\overline{a}x_{i}\overline{y_{i}}\right] }{\left[ 
\func{Re}\left( \overline{a}A\right) \right] ^{\frac{1}{2}}}  \label{5.3} \\
& \leq \frac{1}{2}\cdot \frac{\left\vert A\right\vert +\left\vert
a\right\vert }{\left[ \func{Re}\left( \overline{a}A\right) \right] ^{\frac{1%
}{2}}}\left\vert \sum_{i=0}^{\infty }w_{i}x_{i}\overline{y_{i}}\right\vert .
\notag
\end{align}%
The constant $\frac{1}{2}$ is sharp in (\ref{5.3}).
\end{proposition}

\begin{proof}
Follows by Theorem \ref{t2.1} applied for the inner product $\left\langle
\cdot ,\cdot \right\rangle :\ell _{w}^{2}\left( \mathbb{K}\right) \times
\ell _{w}^{2}\left( \mathbb{K}\right) \rightarrow \mathbb{K}$,%
\begin{equation*}
\left\langle \overline{x},\overline{y}\right\rangle _{w}:=\sum_{i=0}^{\infty
}w_{i}x_{i}\overline{y_{i}}.
\end{equation*}
\end{proof}

\begin{remark}
\label{r5.2}A sufficient condition for (\ref{5.2}) to hold is%
\begin{equation}
\func{Re}\left[ \left( Ay_{i}-x_{i}\right) \left( \overline{x_{i}}-\overline{%
a}\overline{y_{i}}\right) \right] \geq 0\text{ \ for all }i\in \mathbb{N}.
\label{5.4}
\end{equation}
\end{remark}

In the particular case $\rho =1,$ we have the following result.

\begin{corollary}
\label{c5.3}Let $a,A\in \mathbb{K}$ with $\func{Re}\left( \overline{a}%
A\right) >0$ and $\overline{x},\overline{y}\in \ell ^{2}\left( \mathbb{K}%
\right) .$ If%
\begin{equation}
\sum_{i=0}^{\infty }\func{Re}\left[ \left( Ay_{i}-x_{i}\right) \left( 
\overline{x_{i}}-\overline{a}\overline{y_{i}}\right) \right] \geq 0,
\label{5.4b}
\end{equation}%
then one has the inequality%
\begin{align}
\left[ \sum_{i=0}^{\infty }\left\vert x_{i}\right\vert
^{2}\sum_{i=0}^{\infty }\left\vert y_{i}\right\vert ^{2}\right] ^{\frac{1}{2}%
}& \leq \frac{1}{2}\cdot \frac{\sum_{i=0}^{\infty }\func{Re}\left[ A%
\overline{x_{i}}y_{i}+\overline{a}x_{i}\overline{y_{i}}\right] }{\left[ 
\func{Re}\left( \overline{a}A\right) \right] ^{\frac{1}{2}}}  \label{5.5} \\
& \leq \frac{1}{2}\cdot \frac{\left\vert A\right\vert +\left\vert
a\right\vert }{\left[ \func{Re}\left( \overline{a}A\right) \right] ^{\frac{1%
}{2}}}\left\vert \sum_{i=0}^{\infty }x_{i}\overline{y_{i}}\right\vert . 
\notag
\end{align}
\end{corollary}

\begin{remark}
\label{r5.4}If $\mathbb{K}=\mathbb{R}$, then a sufficient condition for
either (\ref{5.1}) or (\ref{5.4}) to hold is%
\begin{equation}
ay_{i}\leq x_{i}\leq Ay_{i}\text{ \ for each \ }i\in \left\{ 1,\dots
,n\right\} ,  \label{5.6}
\end{equation}%
where, in this case, $a,A\in \mathbb{R}$ with $aA>0.$
\end{remark}

When the constants are positive, then the following proposition also holds.

\begin{proposition}
\label{p5.5}Let $m,M>0.$ If $\overline{x},\overline{y}\in \ell
_{w}^{2}\left( \mathbb{K}\right) $ such that 
\begin{equation}
\sum_{i=0}^{\infty }w_{i}\func{Re}\left[ \left( My_{i}-x_{i}\right) \left( 
\overline{x_{i}}-m\overline{y_{i}}\right) \right] \geq 0,  \label{5.7}
\end{equation}%
then one has the inequality%
\begin{equation}
\left[ \sum_{i=0}^{\infty }w_{i}\left\vert x_{i}\right\vert
^{2}\sum_{i=0}^{\infty }w_{i}\left\vert y_{i}\right\vert ^{2}\right] ^{\frac{%
1}{2}}\leq \frac{1}{2}\cdot \frac{M+m}{\sqrt{mM}}\sum_{i=0}^{\infty }w_{i}%
\func{Re}\left( x_{i}\overline{y_{i}}\right) .  \label{5.8}
\end{equation}
\end{proposition}

The proof follows by Corollary \ref{c2.2} applied for the inner product 
\begin{equation*}
\left\langle \overline{x},\overline{y}\right\rangle _{w}:=\sum_{i=0}^{\infty
}w_{i}x_{i},\overline{y_{i}}.
\end{equation*}

The following additive version also holds.

\begin{corollary}
\label{c5.6}With the assumptions in Proposition \ref{p5.5}, one has the
inequalities%
\begin{align}
0& \leq \left[ \sum_{i=0}^{\infty }w_{i}\left\vert x_{i}\right\vert
^{2}\sum_{i=0}^{\infty }w_{i}\left\vert y_{i}\right\vert ^{2}\right] ^{\frac{%
1}{2}}-\sum_{i=0}^{\infty }w_{i}\func{Re}\left( x_{i}\overline{y_{i}}\right)
\label{5.9} \\
& \leq \frac{\left( \sqrt{M}-\sqrt{m}\right) ^{2}}{2\sqrt{mM}}%
\sum_{i=0}^{\infty }w_{i}\func{Re}\left( x_{i}\overline{y_{i}}\right)  \notag
\end{align}%
and%
\begin{align}
0& \leq \sum_{i=0}^{\infty }w_{i}\left\vert x_{i}\right\vert
^{2}\sum_{i=0}^{\infty }w_{i}\left\vert y_{i}\right\vert ^{2}-\left[
\sum_{i=0}^{\infty }w_{i}\func{Re}\left( x_{i}\overline{y_{i}}\right) \right]
^{2}  \label{5.10} \\
& \leq \frac{\left( M-m\right) ^{2}}{4mM}\left[ \sum_{i=0}^{\infty }w_{i}%
\func{Re}\left( x_{i}\overline{y_{i}}\right) \right] ^{2}.  \notag
\end{align}
\end{corollary}

\begin{remark}
\label{r5.7}If $\mathbb{K}=\mathbb{R}$, a sufficient condition for (\ref{5.7}%
) to hold is%
\begin{equation}
my_{i}\leq x_{i}\leq My_{i}\text{ \ for each \ }i\in \mathbb{N},
\label{5.11}
\end{equation}%
where $M>m>0.$
\end{remark}

\end{document}